\documentclass[a4paper, 
11pt
]{article}

\usepackage{subfig}
\usepackage{pgfplots}
\usepackage{pgfplotstable}
\usepackage{mathtools}
\usepackage{multicol}
\usepackage{comment}
\usepackage{booktabs}
\pgfplotsset{compat=1.5}
\usepackage{amssymb}
\usepackage{url}
\usepackage{bm}
\usepackage{symbols}

\usepackage{pdfsync}
\usepackage{float}
\usepackage{tabularx}
\usepackage{enumerate}
\usepackage{array}
\usepackage{xspace}

\usepackage{authblk}

\newcommand{\mupar}{\ensuremath{\boldsymbol{\mu}}}
\newcommand{\etapar}{\ensuremath{\boldsymbol{\eta}}}

\begin{document}

\title{Model Order Reduction by means of Active Subspaces and Dynamic Mode Decomposition for Parametric Hull Shape Design Hydrodynamics}

\author[1]{Marco~Tezzele\footnote{marco.tezzele@sissa.it}}
\author[1,2]{Nicola~Demo\footnote{nicola.demo@sissa.it}}
\author[1]{Mahmoud~Gadalla\footnote{mgadalla@sissa.it}}
\author[1]{Andrea~Mola\footnote{andrea.mola@sissa.it}}
\author[1]{Gianluigi~Rozza\footnote{gianluigi.rozza@sissa.it}}

\affil[1]{Mathematics Area, mathLab, SISSA, International School of Advanced Studies, via Bonomea 265, I-34136 Trieste, Italy}
\affil[2]{Fincantieri - Divisione Navi Mercantili e Passeggeri, Cantieri Navali Italiani SpA, Trieste, Italy}

\maketitle

\begin{abstract}
We present the results of the application of a parameter space
reduction methodology based on active subspaces (AS) to the hull
hydrodynamic design problem. Several parametric deformations of an
initial hull shape are considered to assess the influence of the shape
parameters on the hull wave resistance. Such problem is
relevant at the preliminary stages of the ship design, when
several flow simulations are carried out by the engineers to
establish a certain sensibility with respect to the parameters, which
might result in a high number of time consuming hydrodynamic
simulations. The main idea of this work is to employ the AS to
identify possible lower dimensional structures in the parameter
space. The complete pipeline involves the use of free form deformation
to parametrize and deform the hull shape, the high fidelity solver
based on unsteady potential flow theory with fully nonlinear free
surface treatment directly interfaced with CAD, the use of dynamic
mode decomposition to reconstruct the final steady state given only
few snapshots of the simulation, and the reduction of the parameter
space by AS, and shared subspace. Response surface method is used
to minimize the total drag. 
\end{abstract}

\section{Introduction}
\label{sec:intro}

Nowadays, simulation-based design has naturally evolved into
simulation-based design optimization thanks to new computational
infrastructures and new mathematical methods. In this work we present
an innovative pipeline that combines geometrical parametrization, different model
reduction techniques, and constrained optimization. 
The objective is to minimize the total resistance of a hull advancing in calm water
subject to a constraint on the volume of the hull. We employ free form
deformation (FFD)~\cite{sederbergparry1986,rozza2013free} to parametrize and deform the bottom part of the
stern of the DTMB~5415. For a simulation-based design optimization of
that hull see for example~\cite{serani2016ship}. We select the displacement of some FFD control
points as our parameters, and we sample this parameter space in order
to reduce its dimension by finding an active subspace~\cite{constantine2015active}. In particular we seek a shared
subspace~\cite{ji2018shared} between the target function to minimize
and the constraint function. This subspace allows us to easily perform
the minimization without violating the constraint. As fluid dynamic
model we use a fully nonlinear potential flow one, implemented in
the software WaveBEM
(see~\cite{molaEtAl2013,MolaHeltaiDeSimone2017}). It is interfaced
with CAD data structures, and automatically generates the
computational grids and carry out the simulation with no need for
human intervention. We further accelerate the unsteady flow
simulations through dynamic mode decomposition
presented in~\cite{schmid2010dynamic,kutz2016dynamic} and implemented using
PyDMD~\cite{demo18pydmd}. It allows to reconstruct and predict
all the fields of interest given only few snapshots of the
simulation. 
The particular choice of target function and
constraint does not represent a limitation since the methodology
we present does not rely on those particular
functions. Also the specific part of the domain to be deformed has been
chosen to present the pipeline and does not represent a limitation in
the application of the method.

\section{A Benchmark Problem: Estimation of the Total Resistance}
The hull we consider is
the DTMB~5415, since it is a benchmark for naval hydrodynamics
simulation tools due to the vast experimental data available in the
literature~\cite{olivieri2001towing}. A side view of the
complete hull (used as reference domain $\Omega$) is depicted in Figure~\ref{fig:ffd_hull}.

Given a set of geometrical parameters $\mupar \in \mathbb{D} \subset \mathbb{R}^m$ with
$m \in \mathbb{N}$, we can define a shape morphing function
$\mathcal{M}(\boldsymbol{x}; \boldsymbol{\mu}): \mathbb{R}^3 \to 
\mathbb{R}^3$ that maps the reference domain $\Omega$ to the deformed
one $\Omega(\mupar)$, as $\Omega(\mupar) =
\mathcal{M}(\Omega; \mupar)$.
A detailed description of the specific $\mathcal{M}$ and $\mupar$ used
is in Section~\ref{sec:ffd}.
In the estimation of the total resistance, the simulated flow
field past a surging ship depends on the specific parametric hull shape considered. Thus, the output
of each simulation depends on the parameters defining the deformed shape. To investigate the effect of the shape
parameters on the total drag, we identify a suitable
set of sampling points in $\mathbb{D}$, which, through
the use of free form deformation, define a corresponding set of hull
shapes. Each geometry in such set  is used to run an unsteady fluid dynamic simulation based on a fully nonlinear potential fluid model. 
As a single serial simulation requires approximatively 24h to converge to a steady state solution, 
DMD is employed to reduce such cost to roughly 10h.
The relationship between each point in $\mathbb{D}$ and the estimate for the resistance
is then analyzed by means of AS in order to verify if a further reduction in the parameter space
is feasible. 

\section{Shape Parametrization and Morphing through Free Form Deformation}
\label{sec:ffd}
The free form deformation (FFD) is a widely used technique to deform in a
smooth way a geometry of interest. This section presents a summary of
the method. For a deeper insight on the formulation and more
recent works the reader can refer
to~\cite{sederbergparry1986,lombardi2012numerical,rozza2013free,sieger2015shape,
forti2014efficient,salmoiraghi2016isogeometric,tezzele2017dimension}.
\begin{figure}[htb]
\centering
\includegraphics[width=0.95\textwidth, trim=0 20 0 20]{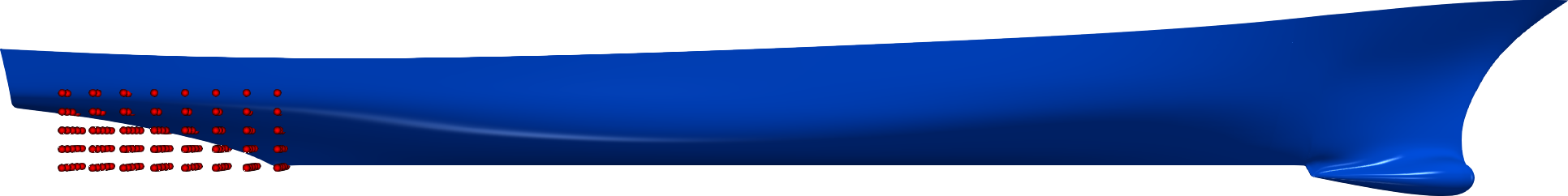}
    \caption{Reference domain $\Omega$, the original DTMB 5415 hull, and the FFD points.}
\label{fig:ffd_hull}
\end{figure}

Basically the FFD needs a lattice of points, called FFD control points, surrounding the object to
morph. Then some of these control points are moved and all the points
of the geometry are deformed by a trivariate tensor-product of
B\'ezier or B-spline functions. The displacements of the control
points are the design parameters $\mupar$ mentioned above. The transformation is composed by 3 different
functions. First we map the physical domain $\Omega$ onto
the reference one $\widehat{\Omega}$ using the map~$\psi$. Then we move some FFD control points $\boldsymbol{P}$ through the map
$\widehat{T}$. This deforms all the points inside the lattice. Finally
we map back to the physical space the reference deformed
domain $\widehat{\Omega}(\mupar)$, obtaining $\Omega(\mupar)$ with $\psi^{-1}$. The composition of these 3 maps is
$T(\cdot, \mupar) = (\psi^{-1} \circ \widehat{T} \circ \psi) (\cdot, \mupar)$.
In Figure~\ref{fig:ffd_hull} we see the lattice of
points around the bottom part of the stern of the
DTMB~5415.
In particular 
we are going to move 7 of them in the vertical direction and 3 along the span
of the boat, so $\mupar \in \mathbb{D} \subset \mathbb{R}^{10}$,
where $\mathbb{D} = [-0.6, 0.5]^{10}$. The original hull corresponds
to $\mupar = 0$. 
We implemented all the algorithms in a Python
package called PyGeM~\cite{pygem}.

\section{High Fidelity Solver based on Fully Nonlinear Potential Flow Model}
\label{sec:solver}
The mathematical model adopted for the simulations is based on potential flow theory. Under the assumptions
of irrotational flow and non viscous fluid, the velocity field admits a scalar potential in
the simply connected flow domain representing the volume of
water surrounding the ship hull. In addition, the Navier-Stokes equation of
fluid mechanics can be simplified to the Laplace equation, which is solved
to evaluate the velocity potential, and to the Bernoulli equation, which allows
the computation of the pressure field. The Laplace equation is
complemented by non penetration boundary condition on the hull, and by fully nonlinear
and unsteady boundary conditions on the water free
surface, written in semi-Lagrangian form~\cite{beck1994}.
In the resulting nonlinear time dependent boundary value problem, the hull is assumed to be
a rigid body moving under the action of gravity and of the hydrodynamic forces
obtained from the pressures resulting from the solution of Bernoulli equation.
The equations governing the hull motions are the 3D
rigid body equations in which the angular displacements are expressed by
unit quaternions.
At each time instant, 
the unknown potential and node displacement fields are obtained by a
nonlinear problem, which results from the spatial and temporal
discretization of the continuous boundary value problem. The spatial discretization of the Laplace
problem is based upon a boundary element method (BEM) described
in~\cite{molaEtAl2013,Giuliani2015}. The domain boundary is
discretized into quadrilateral cells and bilinear shape functions are used to approximate the surface geometry, the velocity
potential values, and the normal component of its surface gradient. The iso-parametric
BEM developed is based on collocating a boundary integral equation~\cite{brebbia}
in correspondence with each node of the computational grid, and on
computing a numerical approximation of the integrals appearing in such
equations. The resulting linear algebraic equations are then combined
with the ODEs derived from the finite element spatial 
discretization of the unsteady fully nonlinear free surface boundary
conditions. 
The final FSI problem is obtained by complementing the 
described system with the equations
of the rigid hull dynamics.
The fully coupled system solution is integrated over time by an arbitrary
order and arbitrary time step implicit backward difference formula scheme.
The potential
flow model is implemented in a C++ software~\cite{molaEtAl2013}. It is equipped with a mesh module directly interfaced
with CAD data structures based on the methodologies for surface mesh generation~\cite{dassi2014}. Thus, for each IGES geometry tested, the
computational grid is generated in a fully automated fashion at the
start of each simulation.
At each time step the wave resistance is computed as $R^w =
\int_{\Gamma^b}p\nb\,d\Gamma\cdot\eb_X$ making use of the pressure
$p$ obtained 
plugging the computed potential in the Bernoulli equation. The non viscous fluid dynamic model drag prediction
is complemented by an estimation of the viscous drag obtained by the
ITTC-57 formula~\cite{MorallITTC1957}.
Results shown in~\cite{MolaHeltaiDesimone2016} indicate that for
Froude numbers in $[0.2, 0.4]$ the total
drag computed for the DTMB~5415 hull differs less than 6\% with
respect to the measurements in~\cite{olivieri2001towing}. For $\text{Fr}=0.28$, at which the present
campaign is carried out, the predicted drag is 46.389~N, which is off by 2.7\% from the correponding 45.175~N experimental
value.  It is reasonable to infer that for
each parametric deformations of the hull the accuracy of the full order model
prediction will be similar to that of the results discussed.

\section{Dynamic Mode Decomposition for Fields Reconstruction}
\label{sec:dmd}
The dynamic mode decomposition (DMD) technique is a tool for the
complex data systems analysis, initially developed in~\cite{schmid2010dynamic} for
the fluid dynamics applications. The DMD provides an approximation of the
Koopman operator capable to describe the system evolution as
linear combination of few linear evolving structures, without requiring any
information about the considered system. We can estimate the future
evolution of these structures in order to reconstruct the system dynamics also
in the future. In this work, we reduce the temporal window where the full order
solutions are computed and we reconstruct the system evolution, applying the
DMD to the output of the full-order model, to gain a significant
reduction of the computational cost.

We define the operator $\mathbf{A}$ such that $x_{k+1} = \mathbf{A} x_{k}$,
   where $x_k$ and $x_{k+1}$ refer respectively to the system state at two
   sequential instants.
To build this operator, we collect several vectors
$\{\mathbf{x}_i\}_{i=1}^m$ that contain the system states equispaced in time,
    called \textit{snapshots}. Let assume all the snapshots have the
    same dimension $n$ and the number of snapshots $m<n$. We can arrange the snapshots in two matrices
$\mathbf{S} = 
\begin{bmatrix}
\mathbf{x}_1 & \mathbf{x}_2 & \dotsc & \mathbf{x}_{m-1}
\end{bmatrix}$ and
$\mathbf{\dot{S}} =
\begin{bmatrix}
\mathbf{x}_2 & \mathbf{x}_3 & \dotsc & \mathbf{x}_{m}
\end{bmatrix}$, 
with $\mathbf{x}_i = \begin{bmatrix} x_i^1 & x_i^2 & \cdots & x_i^n
\end{bmatrix}^\intercal$.
The best-fit $\mathbf{A}$ matrix is given by $\mathbf{A} = \mathbf{\dot{S}}
\mathbf{S}^\dagger$, where $^\dagger$ denotes the Moore-Penrose pseudo-inverse. The biggest issue is
related to the dimension of the snapshots: usually in a complex system the
number of degrees of freedom is high, so the operator $\mathbf{A}$ is very
large. To avoid this, the DMD technique projects the snapshots onto the low-rank
subspace defined by the proper orthogonal decomposition modes.
We decompose the matrix $\mathbf{S}$ using the truncated SVD, that is
$\mathbf{S} \approx \mathbf{U}_r \bm{\Sigma}_r \mathbf{V}^*_r$, and we call
$\mathbf{U}_r$ the matrix whose columns are the first $r$ modes.
Hence, the reduced operator is computed as:
$\mathbf{\tilde{A}} = \mathbf{U}_r^* \mathbf{A} \mathbf{U}_r =
\mathbf{U}_r^* \mathbf{\dot{S}} \mathbf{V}_r \bm{\Sigma}_r^{-1}$. 
We can compute the eigenvectors and eigenvalues of $\mathbf{A}$ through the
eigendecomposition of $\mathbf{\tilde{A}}$, to simulate the system
dynamics. Defining $\mathbf{W}$ and $\bm{\Lambda}$ such that
$\mathbf{\tilde{A}}\mathbf{W} = \bm{\Lambda} \mathbf{W}$, the elements in
$\bm{\Lambda}$ correspond to the nonzero eigenvalues of $\mathbf{A}$ and the
eigenvectors $\bm{\Theta}$ of matrix $\mathbf{A}$, also called DMD
\textit{modes}, can be computed as $\bm{\Theta} = \mathbf{\dot{S}}\mathbf{V}_r
\bm{\Sigma}_r^{-1} \mathbf{W}$.
We implement the algorithm described above, and its most popular
variants, in an open source Python package called PyDMD~\cite{demo18pydmd}. We use
it to reconstruct the evolution of the fluid dynamics system
presented above.

\section{Parameter Space Reduction by means of Active Subspaces}
\label{sec:active}

The active subspaces (AS) property~\cite{constantine2015active} is an
emerging technique for dimension reduction in the parameter
studies. AS has been exploited in several parametrized engineering
models~\cite{grey2017active, constantine2017time, demo2018efficient,
tezzele2017combined}. Considering a multivariate scalar function $f$
depending on the parameters $\mupar$, AS seeks a set of
important directions in the parameter space along which $f$ varies the most. Such directions are linear
combinations of the parameters, and span a lower dimensional
subspace of the input space. This corresponds to a rotation of the
input space that unveils a low dimensional structure of $f$. In the following
we review the AS theory (see~\cite{lukaczyk2014active, constantine2015active}). 
Consider a differentiable, square-integrable scalar function $f
(\mupar): \mathbb{D} \subset \mathbb{R}^m \rightarrow \mathbb{R}$, and a uniform probability density function $\rho: \mathbb{D}
\rightarrow \mathbb{R}^+$. First, we
scale and translate the inputs to be centered at 0 with equal
ranges.
To determine the important directions that most effectively
describe the function variability, the eigenspaces of the
uncentered covariance matrix $\mathbf{C} = \int_{\mathbb{D}} (\nabla_{\mupar} f) ( \nabla_{\mupar} f )^T
\rho \, d \mupar$, needs to be established. 
$\mathbf{C}$ is symmetric positive
definite so it admits a real eigenvalue decomposition, $\mathbf{C} = \mathbf{W} \bm{\Lambda} \mathbf{W}^T,$ 
where $\mathbf{W}$ is a $m \times m$ column matrix of eigenvectors, and
$\bm{\Lambda}$ is the diagonal matrix of non-negative eigenvalues
arranged in descending order. Low eigenvalues suggest that
the corresponding vectors are in the null space of the covariance
matrix, and we can discard those vectors to form an
approximation. The lower dimensional parameter subspace is formed
by the first $M < m$ eigenvectors that correspond to the
relatively large eigenvalues. 
We can partition
$\mathbf{W}$ into $\mathbf{W}_1$ containing the first $M$ eigenvectors which span the
active subspace, and $\mathbf{W}_2$ containing the
eigenvectors spanning the inactive subspace. The dimension reduction
is performed by projecting the full parameter space onto the
active subspace obtaining the active variables
$\mupar_M = \mathbf{W}_1^T\mupar \in \mathbb{R}^M$. The inactive variables are $\etapar = \mathbf{W}_2^T \mupar \in \mathbb{R}^{m - M}$.
Hence, $\mupar \in
\mathbb{D}$ can be expressed as $
\mupar = \mathbf{W}_1\mathbf{W}_1^T\mupar +
\mathbf{W}_2\mathbf{W}_2^T\mupar = \mathbf{W}_1 \mupar_M +
\mathbf{W}_2 \etapar$. The function $f (\mupar)$ can then be
approximated with $g (\mupar_M)$, and the evaluations of some chosen samples $g_i$ for $i = 1, \dots, p \leq N_s$ can be exploited to 
construct a response surface $\mathcal{R}$.

We use the concept of shared subspace~\cite{ji2018shared}. It links the AS of
different functions that share the same parameter space. Expressing
both the objective function and a constraint using the same
reduced variables leads to an easy constrained optimization via the
response surfaces.
The shared subspace $\mathbf{Q}$ between some
$f_i$, $i \in \mathbb{N}$, having an active
subspace of dimension $M$, is defined as follows.
Let us assume that the functions are exactly represented by their AS
approximations, 
then for all $\mathbf{Q} \in \mathbb{R}^{m \times M}$ such that
$\mathbf{W}_{1, \,i}^T \mathbf{Q} = \mathbf{Id}_M$, we have $f_i (\mupar) =f_i (\mathbf{Q}\mathbf{W}_{1, \,i}^T \, \mupar).$
A system of equations needs to be solved for $\mathbf{Q}$, and it can
be proven that $\mathbf{Q}$ will be a linear combination of the active
subspaces of $f_i$.

\section{Numerical Results}
\label{sec:results}

In this section we present the numerical results obtained with the
application of the complete pipeline, presented above, to the DTMB 5415 model hull.

\begin{figure}[htb]
\centering
\includegraphics[width=.27\textwidth, trim=0 0 0 0]{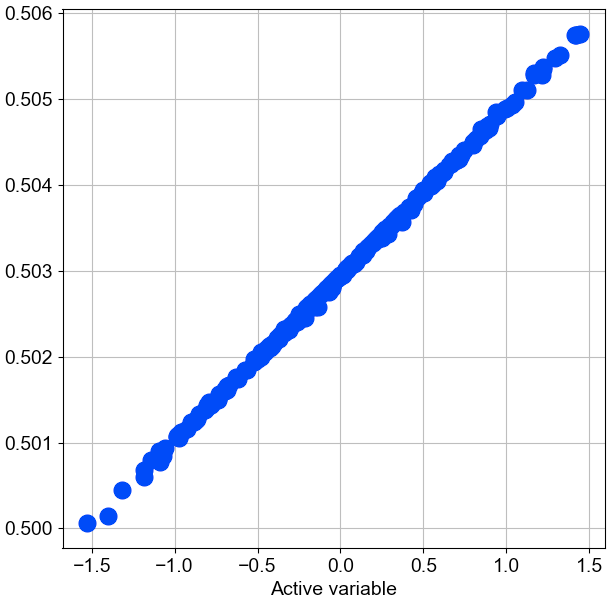}\hfill
\includegraphics[width=.27\textwidth, trim=0 0 0 0]{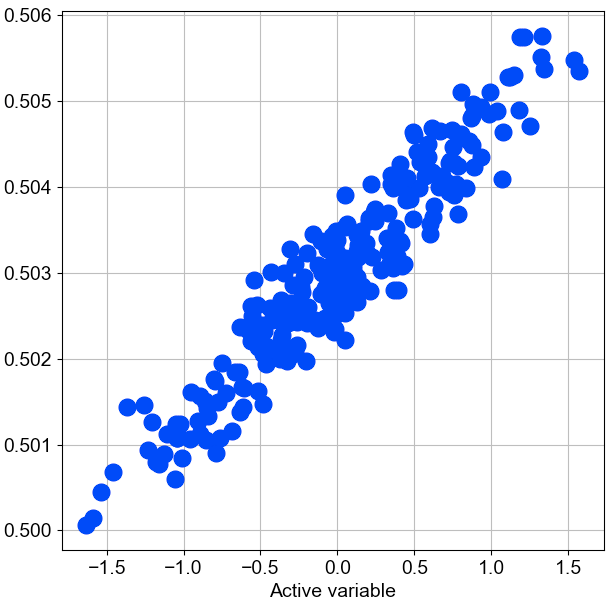}\hfill
\includegraphics[width=.42\textwidth, trim=0 0 0 0]{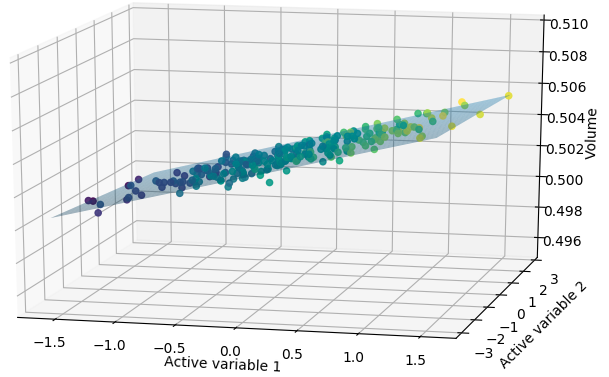}
\caption{On the left the load volume with respect to the active
  variable $\mupar_{\text{vol}}$, while in the center with
  respect to the shared variable $\mupar_Q$. On the right the shared
  variable is 2D and there is a response surface of order one.}
\label{fig:active_shared_vol}
\end{figure}

Using the FFD algorithm implemented in PyGeM~\cite{pygem} we create
200 different deformations of the original hull, sampling uniformly the
parameter space $\mathbb{D} \subset \mathbb{R}^{10}$. 
Each IGES geometry produced is the
input of the full order simulation, in which the hull has been set to advance in calm
water at a constant speed corresponding to Fr $= 0.28$. The full
order computations simulate only 15s of the flow past the
hull after it has been impulsively started from rest. For each
simulation we save the snapshots of the full flow field every 0.1s
between the 7th and the 15th second. The DMD algorithm
implemented in PyDMD~\cite{demo18pydmd} uses these snapshots to
complete the fluid dynamic simulations until convergence to the regime
solution. The reconstructed flow field is then used to calculate the hull
total resistance, that is the quantity of interest we want to
minimize. For each geometry we also compute the volume of the
hull below a certain height $z$ equal for all the hulls. This
is intended as the load volume.
With all the input/output pairs for both the total resistance and the
load volume we can extract the active subspaces for each target
function and compute the shared subspace. Using the shared subspace
has the advantage to allow the representation of the target functions
with respect to the same reduced parameters. The drawback is loosing
the optimality of AS, since it means to shift the
rotation of the parameter space from the optimal one given by
AS. This is clear in Figure~\ref{fig:active_shared_vol},
where the load volume is expressed versus the 1D active variable
$\mupar_{\text{vol}} = \mathbf{W}_{1, \,\text{vol}}^T \,
\mupar$ (on the left), and versus the shared variable
$\mupar_Q = \mathbf{Q}^T\mupar$ in 1D and 2D. The
values of the target function are not perfectly aligned anymore along
the shared variable.

\begin{figure}[htb]
\centering
\includegraphics[width=.32\textwidth, trim=0 10 0 45]{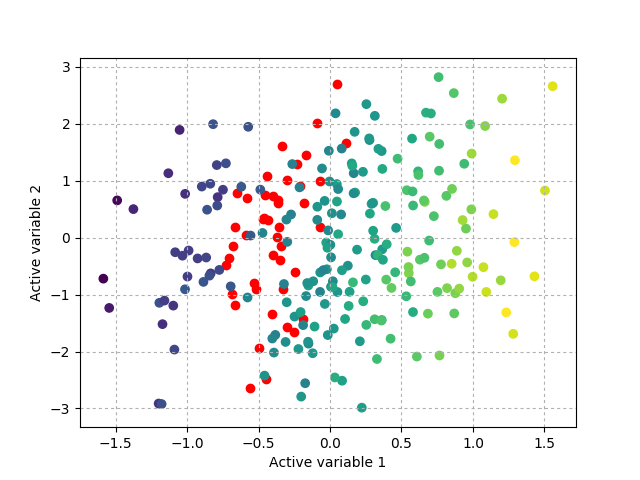}\hfill
\includegraphics[width=.33\textwidth, trim=0 10 0 35]{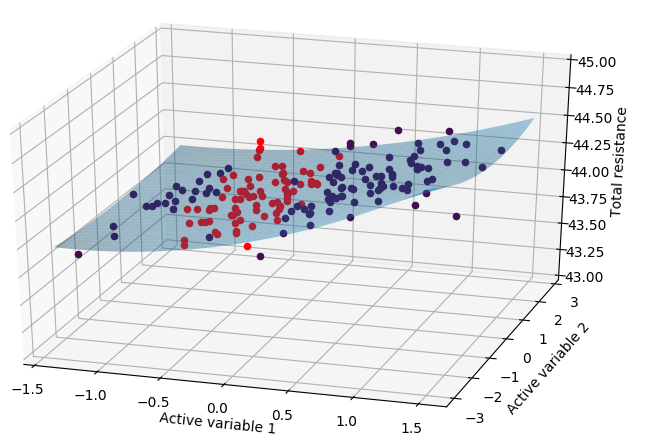}\hfill
\includegraphics[width=.34\textwidth, trim=0 10 0 60]{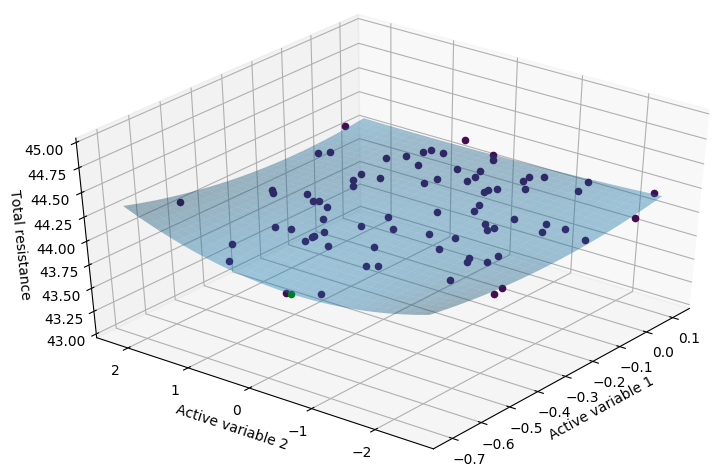}
\caption{On the left and in the center the sufficient summary plot of the resistance
  along with a polynomial response surface. In
  red the points that satisfy the volume constraint. On the right the
  response surface of order two constructed with the red points. In
  green the minimum.}
\label{fig:resistance_surf}
\end{figure}

To capture the most information while having the possibility
to plot the results, we select the active subspace to be of dimension
2 for both quantities of interest and we compute the
corresponding shared subspace, that is a linear combination of the
active subspaces of the two functions.
Then we select a lower and an upper bound in which we constraint the
load volume. After, we compute the subset of the shared subspace that satisfies
such constraint in order to impose it on the total resistance. In
Figure~\ref{fig:resistance_surf} we plot on the left the volume against the
bidimensional shared variable (seen from above) and in red the
realizations of the volume between the imposed bounds. 
In the center we plot the sufficient summary plot for the total resistance
with respect to the shared variables and highlight the simulations
that satisfy the volume constraint in red, together with the response
surface. We notice that here the data are more
scattered if compared with the volume. This is due to the high
nonlinearity of the problem and the inactive directions we
discarded. On the right we construct a polynomial response
surface of order 2 using only the red dots and in green we
highlight the minimum of such response surface. That represents the
minimum of the total resistance in the reduced parameter space
subjected to the volume constraint. The root mean square error
committed with the response surface is around 0.7 Newton
depending on the run. We can map such reduced point in the full space of the parameters and
deform the original hull accordingly. This represents a good starting
point for a more sophisticated optimization. 

\section{Conclusions and Perspectives}
\label{sec:the_end}

In this work we presented a complete pipeline composed of shape
parametrization, hydrodynamic simulations, and model reduction
combining DMD and AS. We applied it to the problem of minimizing the
total resistance of a hull advancing in calm water, subject to a load
volume constraint, varying the shape of the DTMB~5415. We expressed the
two functions in the same reduced space and we constructed a response
surface to find the minimum of the total drag that satisfies the
volume constraint. We committed an error around 0.7 Newton with respect to
the full order solver. This minimum can be used as a starting
point of a more sophisticated optimization algorithm. The proposed pipeline is
independent of the specific deformed part of the hull, of the
solver used, and of the function to minimize. 
Future work will focus on several different areas to improve physical and
mathematical aspects of the algorithms presented, as well as to better integrate its
different parts, and to automate the simulation campaign with post
processing processes.

\section*{Acknowledgements}
This work was partially supported by the project PRELICA,
``Advanced methodologies for hydro-acoustic design of naval propulsion'', supported by Regione
FVG, POR-FESR 2014-2020, Piano Operativo Regionale Fondo Europeo per lo Sviluppo Regionale,
partially funded by the project HEaD, ``Higher
Education and Development'', supported by Regione FVG, European
Social Fund FSE 2014-2020, and by European Union
Funding for Research and Innovation --- Horizon 2020 Program --- in
the framework of European Research Council Executive Agency: H2020 ERC
CoG 2015 AROMA-CFD project 681447 ``Advanced Reduced Order Methods
with Applications in Computational Fluid Dynamics'' P.I. Gianluigi
Rozza.


\end{document}